\documentclass[twocolumn,twoside]{IEEEtran}
\usepackage{amsmath,amssymb,amsthm,wasysym,epsfig,dsfont,color,subfigure,empheq,graphicx}
\usepackage{enumerate,url,algorithm,empheq,balance,url,multirow}
\usepackage{comment}
\usepackage{algpseudocode}
\usepackage[dvipsnames]{xcolor}
\usepackage[utf8]{inputenc}
\usepackage{cite}
\usepackage{threeparttable}
\usepackage{ifthen}
\usepackage{hyperref}
\usepackage{tikz} 
\usepackage{pgfplots}
\usetikzlibrary{arrows}
\usetikzlibrary{plotmarks}
\hypersetup{colorlinks=true,linkcolor=black,urlcolor=black}
\usepackage{booktabs}
\usepackage{cases}
\usepackage{enumitem}

\usepackage[numbers,sort&compress]{natbib}

\usepackage{placeins}
\usetikzlibrary{shadows} % LaTeX and plain TeX
\usetikzlibrary[shadows] % ConTeXt
\usetikzlibrary{patterns}
\usetikzlibrary{calc}

\hypersetup{citecolor=black}

\definecolor{Eggplant}{rgb}{0.235, 0.106, 0.557}

\DeclareMathOperator*{\Xmin}{Xmin}
\DeclareMathOperator*{\Xmax}{Xmax}
\DeclareMathOperator*{\Qmin}{Qmin}
\DeclareMathOperator*{\Qmax}{Qmax}

\DeclareMathOperator*{\UP}{UP}
\DeclareMathOperator*{\DN}{DN}

\allowdisplaybreaks

\newenvironment{VarDescription}[1]%
{\begin{list}{}{%
			\settowidth{\labelwidth}{\textbf{#1}}%
			\setlength{\leftmargin}{\labelwidth}\addtolength{\leftmargin}{\labelsep}}}%
	{\end{list}}

\begin{document}
	
	\title {Impact of Ramping Constraints on Electricity Market Outcomes}
	
	\author{
		Xuan Liu, \IEEEmembership{Member,~IEEE,} and 
		Antonio\ J.\ Conejo, \IEEEmembership{Fellow,~IEEE}
		\vspace*{-20pt}
		\par		
	}
	
	\markboth{TBD}{Liu, X., Conejo, A.~J.}
	
	\maketitle
	
	\begin{abstract}
		High penetration of renewable energy brings the need of extra ramping capability. It is thus critical to understand the impact of ramping constraints on electricity market outcomes, and particularly on prices. Specifically, to study the market clearing price formation with ramping constraints, we propose a model that can identify meaningful market equilibria. We assume that producers/consumers offer/bid strategically both quantities and prices. In turn, the market operator clears the market with collected offers, bids, and ramping constraints of producers, which results in market outcomes, including prices. We use an illustrative example and a case study to analyze the impact of ramping constraints on market outcomes in competitive and strategic cases.
	\end{abstract}
	
	
	\section*{Notation}
	\label{Notation}
	%\vspace{-3.5 mm}
	The notation used in this paper is provided below for easy reference.
	
	\noindent \textit{A. Indices and Sets:} 
	\begin{VarDescription}{$\hspace{10mm}$}\setlength\itemsep{1mm}
		\item [{$i$}] {Index for producers running from 1 to $I$.}
		\item [{$j$}] {Index for consumers running from 1 to $J$.}
		\item [{$t$}] {Index for time periods running from 1 to $T$}
	\end{VarDescription}
	%\vspace{+1 mm}
	\noindent \textit{B. Parameters:}
	\begin{VarDescription}{$\hspace{10mm}$}\setlength\itemsep{1mm}
		\item [{$c_i$}] {Marginal cost of producer $i$.}
		\item [{$u_{j,t}$}] {Marginal utility of consumer $j$ in period $t$.}
		\item [{$c^{\rm H}$}] {Upper bound of offer prices for producers.}
		\item [{$u^{\rm L}$}] {Lower bound of bid prices for consumers.}
		\item [{$X_i^{\min}$}] {Minimum production level of producer $i$.} 
		\item [{$X_{i}^{\max}$}] {Maximum production level of producer $i$.} 
		\item [{$Q_j^{\min}$}] {Minimum consumption level of consumer $j$.} 
		\item [{$Q_{j,t}^{\max}$}] {Maximum consumption level of consumer $j$ in period $t$.}
		\item [{$R_i^{\UP}$}] {Ramping-up limit of producer $i$.}
		\item [{$R_i^{\DN }$}] {Ramping-down limit of producer $i$.}
		\item [{$x_{i,0}$}] {Initial production level of producer $i$.}
		\item [{$\widetilde{X}_{i,0}$}] {Initial production cap of producer $i$.}
        \item [{$\rho$}] {A small enough scalar, such as 0.01.}
	\end{VarDescription}
	\noindent \textit{C. Variables:}
	\begin{VarDescription}{$\hspace{10mm}$}\setlength\itemsep{1mm}
		\item [{$x_{i,t}$}] {Production of producer $i$ in period $t$.}
		\item [{$q_{j,t}$}] {Consumption of consumer $j$ in period $t$.}
		\item [{$c_{i,t}^{\rm S}$}] {Strategic offer price of producer $i$ in period $t$.}
		\item [{$u_{j,t}^{\rm S}$}] {Strategic bid price of consumer $j$ in period $t$.}
		\item [{$\widetilde{X}_{i,t}$}] {Strategic production cap of producer $i$ in period $t$.}
		\item [{$\widetilde{Q}_{j,t}$}] {Strategic consumption cap of consumer $j$ in period $t$.}
		\item [{$ \lambda_{t} $}] {Locational marginal price (LMP) in period $t$.}
	\end{VarDescription}

	\noindent \textit{D. Acronyms:}
	\begin{VarDescription}{$\hspace{10mm}$}\setlength\itemsep{1mm}
		\item [{RC}] {Ramping Constraint.}
		\item [{MO}] {Market Operator.}
		\item [{LMP}] {Locational Marginal Price.}
		\item [{KKT}] {Karush–Kuhn–Tucker.}
		\item [{FRP}] {Flexible Ramping Product.}
		\item [{SW}] {Social Welfare.}
		\item [{TP}] {Total Profit.}
		\item [{TS}] {Total Surplus.}
	\end{VarDescription}

	\section{Introduction}
	\label{Introduction}

	The increasing penetration of renewable resources (solar and wind power units) brings clean energy but also challenges the operation of power systems. A high proportion of renewable energy increases the slope of change of the net-load curve \cite{Denholm(2015)}. This curve reveals the need for flexible dispatchable resources, which can ramp their output up/down quickly, to ensure power balance \cite{Conejo(2018)}. Thus, it is necessary to investigate the impact of ramping constraints (RCs) on the market outcome and specially on the price formation. This analysis may provide a reference for future market design.
	
	In this paper, we propose a multi-period equilibrium model to identify market equilibria that extends and enhances the model proposed in \cite{Liu(2021)}.
	
	We consider a day-ahead energy market cleared by a market operator (MO) once a day, and on an hourly basis. The market clearing outcome provides hourly quantities to be produced and consumed, and the locational marginal prices (LMPs). Each strategic producer/consumer submits offer/bid prices and quantity caps seeking its maximum profit/surplus. 
	The offer/bid behavior of each producer/consumer is modeled as an optimization problem. 
	The MO clears the market to maximize the social welfare based on the collected offers and bids, which results in an optimization problem as well.
	Producers' RCs are explicitly included in the MO problem. With interrelated decision variables, all these optimization problems compose an equilibrium problem \cite{Conejo(2020)}.  
	
	Since the considered optimization problems are linear, they can be replaced by their corresponding Karush-Kuhn-Tucker (KKT) optimality conditions, which are linear and nonlinear constraints. 
	These constraints can be jointly solved using an auxiliary optimization problem with pre-specified objective functions (OF) to obtain meaningful equilibria. 
	
	\subsection{Literature Review and Contributions}	
	The existing literature analyze the impact of RCs mainly based on flexible ramping products (FRPs). 
	Reference \cite{Zou(2016a)} and \cite{Chen(2017)} use Nash-Cournot equilibrium models to assess how FRPs affect the market equilibrium and concludes that energy prices increase slightly with FRPs. 
	Reference \cite{Cornelius(2014)} proposes unit commitment and economic dispatch models to evaluate FRPs impact from the MO's viewpoint. 
	It shows that FRPs can reduce price spikes caused by flexibility shortage during ramping intervals, which offset raised prices in non-shortage intervals, and may lead to lower average prices overall. 
	Reference \cite{Cornelius(2018)} points out that the performance of FRPs is highly sensitive to parameters setting and is very system dependent. 
	Reference \cite{Wang(2014)} proposes deterministic and stochastic models to evaluate the effect of FRPs and concludes that FRPs should be carefully designed to ensure that they can effectively improve the market efficiency. 
	
	Some other works investigate strategic behaviors of specific agents related to ramping limits. 
	An equilibrium problem with equilibrium constraints (EPEC) is proposed in \cite{Moiseeva(2015)} to investigate fast-ramping generators' gaming behaviors regarding ramp rates. 
	The authors conclude that strategic generators are prone to exercise market power by offering lower ramp rates in cases that require flexibility. 
	Two types of market structures are compared in \cite{Moiseeva(2017)} to study which one of them is better to incentivize generators to offer true ramp rates. 
	The authors conclude that in a duopoly case, the ramp rate offer in a one-stage market (ramp rates and productions are offered simultaneously) is generally higher than that in a two-stage market (ramp rates and productions are offered sequentially). 
	A mathematical program with equilibrium constraints (MPEC) is proposed in \cite{Nasrolahpour(2018)} to study the impact of RCs on the strategic behavior of energy storage facilities. 
	The authors conclude that strategic storage facilities benefit from the supply-side inflexibility cause by RCs, especially during peak hours. 
	Reference \cite{Shahmohammadi(2018)} examines the effectiveness of mitigating ramping effects by energy storage with different types of market participation. 
	The authors conclude that co-ownership and co-operation of renewable generators and energy storage is the best way to mitigate market inefficiencies caused by RCs. 
	
	Within the context above, we proposed a multi-period equilibrium model to study the impact of RCs on electricity market outcomes, and on clearing prices in particular. 
	An illustrative example and realistic case studies are analyzed in this paper. 
	The outcomes of these cases allow comprehending the impact of RCs on the market clearing price, with and without market power.
	
	The main contribution of this paper is providing a compact and efficient multi-period equilibrium model that allows comprehending the impact of RCs on prices.

	\subsection{Paper Structure}
	
	The rest of this paper is organized as follows. 
	Section \ref{Formulation} provides the equilibrium model formulation and its solution using an auxiliary optimization problem.
	Section \ref{IllustrativeExample} and \ref{CaseStudy} provide results for an illustrative example and a realistic case study, respectively, which are discussed and analyzed in detail. 
	Section \ref{Conclusion} draws some conclusions.

	\section{Formulation}
	\label{Formulation}

	\subsection{Model Assumption}
	\label{Model Assumption}
	
	The main assumptions of the proposed multi-period equilibrium model are summarized below:
	
	\begin{enumerate}
		\item Multiple time periods are considered, but the transmission network is neglected. 	
		This way the impact of network constraints and of intertemporal constraints on the locational marginal price are decoupled. 
		\item A strategic producer/consumer may submit different offer/bid prices and different production/consumption caps in different periods to maximize its own total profit/surplus. 
		\item The market operator maximizes the ``declared'' social welfare for all periods based on offers and bids.
		\item In different periods, the true cost of a producer does not change but the true utility of a consumer does. For simplicity, we consider that the upper/lower bounds of offer/bid prices are the same for all periods.
		\item The locational marginal price of each period is the dual variable of the corresponding power balance constraint.
		\item RCs are included in producers' problems and in the market operator's problem.
		\item Only producers with RCs are considered. 
		\item Consumers do not exhibit RCs.
	\end{enumerate}
	
	\subsection{Multi-Period Single-Level Models}
	\label{Single-Level Model}
	
	The models below are formulated for the three types of participants in the market: producers, consumers, and the market operator. Each producer/consumer submits strategic offer/bid prices and production/consumption caps to the market operator. Then, the market operator clears the market based on the received offers and bids.

	The problem of strategic producer $i$ is:	
	\begin{equation}\tag{1a}\label{eq:1a}
		\underset{\widetilde{X}_{i,t}, c_{i,t}^{\rm S}}{\text{Minimize}}\; \sum_{t}\left(c_i-\lambda_{t}\right)\widetilde{X}_{i,t}
	\end{equation}
	\vspace{-5mm}
	\begin{align}
		\quad &\text{subject to:}&\notag\\
		&c_i \leq c_{i,t}^{\rm S} \leq c^{\rm H} :\mu_{i,t}^{\rm C}, \mu_{i,t}^{\rm  H}  \;\; \forall t\tag{1b}\label{eq:1b}\\
		&X^{\min}_i\leq \widetilde{X}_{i,t} \leq X^{\max}_{i} :\mu_{i,t}^{\Xmin}, 
		\mu_{i,t}^{\Xmax}  \;\; \forall t \tag{1c}\label{eq:1c}\\
		&-R_i^{\DN}\leq \widetilde{X}_{i,t}-\widetilde{X}_{i,0}\leq R_i^{\UP}:\mu_{i,t}^{\DN}, \mu_{i,t}^{\UP} \;  t=1\tag{1d}\label{eq:1d}\\
		&-R_i^{\DN}\leq \widetilde{X}_{i,t}-\widetilde{X}_{i,t-1}\leq R_i^{\UP}: \mu_{i,t}^{\DN}, \mu_{i,t}^{\UP}\;  \forall t> 1 \tag{1e}\label{eq:1e}
	\end{align}	
	where $ \widetilde{X}_{i,t}$ and $c_{i,t}^{\rm S} $ are the primal optimization variables, which represent strategic offers for all periods. 
	
	The objective function (\ref{eq:1a}) is the negative total profit of producer $i$ over all periods. 
	Constraints (\ref{eq:1b}) and (\ref{eq:1c}) stablish ranges of offer prices and production quantity caps, respectively. 
	Constraints (\ref{eq:1d}) and (\ref{eq:1e}) are ramping limits for producer $i$. 
	Note that constraints (\ref{eq:1d}) and (\ref{eq:1e}) need to be include into the producer's problem; otherwise, the intertemporal opportunity cost \cite{Zhao(2020)} cannot be reflected properly. 
	Without constraints (\ref{eq:1d}) and (\ref{eq:1e}), the strategic offers of the producer in each period are independent.

	The problem of strategic consumer $j$ is:	
	\begin{equation}\tag{2a}\label{eq:2a}
		\underset{\widetilde{Q}_{j,t}, u_{j,t}^{\rm S}}{\text{Minimize}}\;  \sum_{t}\left(\lambda_{t}-u_{j,t}\right)\widetilde{Q}_{j,t}
	\end{equation}
	\vspace{-5mm}
	\begin{align}
		\quad &\text{subject to:}&\notag\\
		&u^{\rm L}\leq u_{j,t}^{\rm S} \leq u_{j,t} :\mu_{j,t}^{\rm L}, \mu_{j,t}^{\rm U} \;\; \forall t\tag{2b}\label{eq:2b}\\
		&Q^{\min}_j\leq \widetilde{Q}_{j,t} \leq Q^{\max}_{j,t} :\mu_{j,t}^{\Qmin}, \mu_{j,t}^{\Qmax}  \;\; \forall t \tag{2c}\label{eq:2c}
	\end{align}
	where $ \widetilde{Q}_{j,t}$ and $u_{j,t}^{\rm S} $ are the primal optimization variables, which represent strategic bids for all periods.	The objective function (\ref{eq:2a}) is the negative total surplus of consumer $j$.  Constraints (\ref{eq:2b}) and (\ref{eq:2c}) stablish ranges of bid prices and consumption quantity caps, respectively.
	
	The market operator's problem is:	
	\begin{equation}\tag{3a}\label{eq:3a}
		\underset{x_{i,t}, q_{j,t}}{\text{Minimize}}\; \sum_{i,t} c_{i,t}^{\rm S}x_{i,t}-\sum_{j,t} u_{j,t}^{\rm S}q_{j,t}
	\end{equation}
	\vspace{-5mm}
	\begin{align}
		&\text{subject to:}&\notag\\
		&\sum_{j}q_{j,t} - \sum_{i}x_{i,t} = 0: \lambda_{t}\;\; \forall t\tag{3b}\label{eq:3b}\\
		&X^{\min}_i \leq x_{i,t} \leq \widetilde{X}_{i,t}  :\delta^{\Xmin}_{i,t}, \delta^{\Xmax}_{i,t}\;\;  \forall i, \forall t \tag{3c}\label{eq:3c}\\
		& Q^{\min}_j \leq q_{j,t} \leq \widetilde{Q}_{j,t}   : \delta^{\Qmin}_{j,t}, \delta^{\Qmax}_{j,t} \;\; \forall j, \forall t\tag{3d}\label{eq:3d}\\
		&-R_i^{\DN}\leq x_{i,t}-x_{i,0}\leq R_i^{\UP}:\delta_{i,t}^{\DN}, \delta_{i,t}^{\UP} \; \forall i, t=1\tag{3e}\label{eq:3e}\\
		&-R_i^{\DN}\leq x_{i,t}-x_{i,t-1}\leq R_i^{\UP}: \delta_{i,t}^{\DN}, \delta_{i,t}^{\UP}\; \forall i, \forall t> 1 \tag{3f}\label{eq:3f}		
	\end{align}	
	where $ x_{i,t}$ and $q_{j,t}$ are the primal optimization variables, which represent production and consumption quantities, respectively. 
	
	The objective function (\ref{eq:3a}) is the declared negative social welfare, which is computed based on the submitted offers and bids. 
	The power balance per period is equation (\ref{eq:3b}). 
	Constraints (\ref{eq:3c}) and (\ref{eq:3d}) limit the production and consumption of each producer and consumer, respectively. 
	Constraints (\ref{eq:3e}) and (\ref{eq:3f}) are the ramping limits for producers with RCs.

	\subsection{KKT Conditions}
	\label{KKT Formulation}	
	 
	 Solving jointly interrelated problems (1), (2), and (3) requires that each problem is replaced by its optimality conditions. 

	The KKT conditions of producers' problems are:
	\begin{align}
		&\Big\{(1{\rm b})-(1{\rm e})\Big.\tag{4a}\label{eq:4a}\\
		&c_i-\lambda_{t}-\mu_{i,t}^{\Xmin}+\mu_{i,t}^{\Xmax} +\mu_{i,t}^{\UP}-\mu_{i,t}^{\DN} = 0  \;\; t=T \tag{4b}\label{eq:4b} \\
		&c_i-\lambda_{t}-\mu_{i,t}^{\Xmin}+\mu_{i,t}^{\Xmax}\notag\\
		&\qquad +\mu_{i,t}^{\UP}-\mu_{i,t}^{\DN}-\mu_{i,t+1}^{\UP}+\mu_{i,t+1}^{\DN} = 0  \;\; \forall t<T \tag{4c}\label{eq:4c}\\			
		&\mu_{i,t}^{\rm H}-\mu_{i,t}^{\rm C}= 0  \tag{4d}\label{eq:4d}\\			
		&0\leq \widetilde{X}_{i,t}-X_i^{\min}\perp\mu_{i,t}^{\Xmin} \geq 0 \;\;  \forall t \tag{4e}\label{eq:4e}\\						
		&0\leq -\widetilde{X}_{i,t}+X_i^{\max}\perp\mu_{i,t}^{\Xmax}\geq 0  \;\;  \forall t \tag{4f}\label{eq:4f}\\
		&0\leq c_{i,t}^{\rm S}-c_i \perp \mu_{i,t}^{\rm C}\geq 0  \;\;  \forall t \tag{4g}\label{eq:4g}\\
		& 0\leq -c_{i,t}^{\rm S}+c^{\rm H} \perp \mu_{i,t}^{\rm H}\geq 0 \;\;  \forall t \tag{4h}\label{eq:4h}\\		
		&0\leq \widetilde{X}_{i,t}-\widetilde{X}_{i,0}+R_{i}^{\DN} \perp \mu_{i,t}^{\DN}\geq 0  \;\;  t=1 \tag{4i}\label{eq:4i}\\
		&0\leq -\widetilde{X}_{i,t}+\widetilde{X}_{i,0}+R_{i}^{\UP} \perp \mu_{i,t}^{\UP}\geq 0  \;\;  t=1 \tag{4j}\label{eq:4j}\\
		&0\leq \widetilde{X}_{i,t}-\widetilde{X}_{i,t-1}+R_{i}^{\DN} \perp \mu_{i,t}^{\DN}\geq 0  \;\; \forall t>1 \tag{4k}\label{eq:4k}\\
		&\Big. 0\leq -\widetilde{X}_{i,t-1}+x_{i,0}+R_{i}^{\UP} \perp \mu_{i,t}^{\UP}\geq 0  \;\;  \forall t>1 \Big\}\; \forall i \tag{4l}\label{eq:4l}		
	\end{align}	
	Constraints (\ref{eq:4a})-(\ref{eq:4d}) are primal and dual constraints, and constraints (\ref{eq:4e})-(\ref{eq:4l}) are complementarity conditions.

	The KKT conditions of the consumers' problems are:
	\begin{align}
		&\Big\{(2{\rm b})-(2{\rm c})\Big.\tag{5a}\label{eq:5a}\\
		& \Big\{\Big. \lambda_{t}-u_{j,t}-\mu_{j,t}^{\Qmin}+\mu_{j,t}^{\Qmax}=0  \tag{5b}\label{eq:5b}\\
		& \mu_{j,t}^{\rm U}-\mu_{j,t}^{\rm L} = 0   \tag{5c}\label{eq:5c}\\
		&0\leq \widetilde{Q}_{j,t}-Q_j^{\min} \perp \mu_{j,t}^{\Qmin} \geq 0  \tag{5d}\label{eq:5d}\\
		&0\leq -\widetilde{Q}_{j,t}+Q_{j,t}^{\max} \perp \mu_{j,t}^{\Qmax} \geq 0  \tag{5e}\label{eq:5e}\\
		&0\leq u_{j,t}^{\rm S}-u^{\rm L} \perp \mu_{j,t}^{\rm L}  \geq 0 \tag{5f}\label{eq:5f}\\
		&\Big. 0\leq -u_{j,t}^{\rm S}+u_{j,t} \perp \mu_{j,t}^{\rm U}  \geq 0  \tag{5g}\label{eq:5g}  \Big\}\; \forall t\;  \Big\}\; \forall j 		
	\end{align}		
	Constraints (\ref{eq:5a})-(\ref{eq:5c}) are primal and dual constraints, and constraints (\ref{eq:5d})-(\ref{eq:5g}) are complementarity conditions.

	The KKT conditions of the market operator's problem are:
	\begin{align}
		& (1.3 {\rm b}) - (1.3 {\rm h}) \tag{6a}\label{eq:6a}\\	
		&c_{i,t}^{\rm S}-\lambda_{t}-\delta_{i,t}^{\Xmin}+\delta_{i,t}^{\Xmax} \notag\\
		&\qquad\qquad\qquad\qquad\quad+\delta_{i,t}^{\UP}-\delta_{i,t}^{\DN} = 0  \;\;  \forall i,t=T \tag{6b}\label{eq:6b} \\
		&c^{\rm S}_{i,t}-\lambda_{t}-\delta_{i,t}^{\Xmin}+\delta_{i,t}^{\Xmax}\notag\\
		&\qquad + \delta_{i,t}^{\UP}-\delta_{i,t}^{\DN}-\delta_{i,t+1}^{\UP}+\delta_{i,t+1}^{\DN} = 0  \;\; \forall i,\; \forall t<T \tag{6c}\label{eq:6c}\\
		&\lambda_{t}-u_{j,t}^{\rm S}-\delta_{j,t}^{\Qmin}+\delta_{j,t}^{\Qmax} = 0 \;\;  \forall j,\forall t \tag{6d}\label{eq:6d}\\		
		&0\leq x_{i,t}-X_{i}^{\min}  \perp \delta_{i,t}^{\Xmin} \geq 0 \;\; \forall i,\; \forall  t       \tag{6e}\label{eq:6e}\\
		&0\leq -x_{i,t}+\widetilde{X}_{i,t} \perp \delta_{i,t}^{\Xmax} \geq 0 \;\; \forall i,\; \forall t\tag{6f}\label{eq:6f}\\
		&0\leq q_{j,t}-Q_{j}^{\min}  \perp \delta_{j,t}^{\Qmin} \geq 0 \;\; \forall j,\; \forall t\tag{6g}\label{eq:6g}\\
		&0\leq -q_{j,t}+\widetilde{Q}_{j,t} \perp \delta_{j,t}^{\Qmax} \geq 0 \;\; \forall j,\; \forall t\tag{6h}\label{eq:6h}\\	
		&0\leq x_{i,t}-x_{i,0}+R_i^{\DN} \perp \delta_{i,t}^{\DN} \geq 0  \;\; \forall i,\; t=1 \tag{6i}\label{eq:6i}	\\
		&0\leq -x_{i,t}+x_{i,0}+R_i^{\UP} \perp \delta_{i,t}^{\UP} \geq 0 \;\; \forall i,\; t=1 \tag{6j}\label{eq:6j}	\\
		&0\leq x_{i,t}-x_{i,t-1}+R_i^{\DN} \perp \delta_{i,t}^{\DN} \geq 0  \;\; \forall i,\; \forall t\geq 1 \tag{6k}\label{eq:6k}	\\
		&0\leq -x_{i,t}+x_{i,t-1}+R_i^{\UP} \perp \delta_{i,t}^{\UP} \geq 0\;  \forall i,\forall t\geq 1 \tag{6l}\label{eq:6l} 
	\end{align}
	Constraints (\ref{eq:6a})-(\ref{eq:6d}) are primal and dual constraints, and constraints (\ref{eq:6e})-(\ref{eq:6l}) are complementarity conditions.
	\subsection{Auxiliary Problem}
	\label{Auxiliary Problem}
	
	We use an auxiliary problem to solve all KKT conditions together ((4), (5), and (6)). This problem has the form:
	\begin{equation}\tag{10a}\label{eq:10a}
		\underset{\Xi}{\text{Maximize}}\;\; \rm OF
	\end{equation}
	\vspace{-5mm}
	\begin{align}
		\quad&\text{subject to:}&\notag\\
		&\qquad\quad(4)-(6), \tag{10b}\label{eq:10b}
	\end{align}		
	where set $\Xi$ contains all primal and dual variables of optimization problems (1)-(3). The objective function is selected pursuing different targets. Typical objective functions can be:
	\begin{align}
	&{\rm SW} = \sum_{j,t}u_{j,t}q_{j,t}-\sum_{i,t}c_ix_{i,t} \tag{11a}\label{eq:11a}\\
	&{\rm TP} = \sum_{i,t}\left(\lambda_{t}-c_i\right)x_{i,t} \tag{11b}\label{eq:11b}\\
	&{\rm TS} = \sum_{j,t}\left(u_{j,t}-\lambda_{t}\right)q_{j,t} \tag{11c}\label{eq:11c}\\
	&{\text{IP}(i)} = \sum_{t}\left(\lambda_{t}-c_i\right)x_{i,t} - \rho\sum_{t}c_{i,t}^{\rm S} \tag{11d}\label{eq:11d}\\
	&{\text{IS}(j)} = \sum_{t}\left(u_{j,t}-\lambda_{t}\right)q_{j,t} + \rho\sum_{t}u_{j,t}^{\rm S} \tag{11e}\label{eq:11e}
	\end{align}
	
	Objective function \eqref{eq:11a} is true social welfare (SW). 
	Objective functions \eqref{eq:11b} and \eqref{eq:11c} are the total profit (TP) of producers and the total surplus (TS) of consumers, respectively. 
	Using objective functions \eqref{eq:11b} or \eqref{eq:11c}, we can identify the equilibria that correspond to producers or consumers exercising market power, respectively.
	Objective functions \eqref{eq:11d} and \eqref{eq:11e} are related to individual profits/surplus.
	Using objective function \eqref{eq:11d}, we can identify the periods in which a strategic producer offers a higher price for maximizing its profits.
	Similarly, using objective function \eqref{eq:11e}, we can identify the periods in which a strategic consumer bids a lower price to maximize its surplus.

	The SW (\ref{eq:11a}) is a linear function. On the other hand, the TP, the TS, the $\text{IP}(i)$, and the  $\text{IS}(j)$ are nonlinear functions.
	In addition, we note that the complementarity condition $0\leq x\perp y\geq 0$ is equivalent to $x\geq0$, $y\geq0$, and $xy=0$ \cite{Ben-TalNemirovski(2001)}, \cite{Gabriel(2013)}. 
	Thus, (10) is a continuous nonlinear optimization problem and can be solved by global nonlinear solvers. 
	
		\section{Illustrative Example}
		\label{IllustrativeExample}

	%%%%%%%%%%%%%%%%%%%%%%%%%%%%%%%%%%%%%%%%%%
	%Drawing 
	\tikzset{Solar/.pic={
			code={
				\tikzstyle{every node}=[font=\LARGE]
				\draw [BurntOrange] (0,0) circle (0.5);
				\node [BurntOrange] at (0,0) {S};
				\draw [BurntOrange] (0,-0.5) -- (0,-1);
			}
		}
	}
	
	\tikzset{Thermal/.pic={
			code={
				\tikzstyle{every node}=[font=\LARGE]
				\draw  (0,0) circle (0.5);
				\node  at (0,0) {T};
				\draw  (0,-0.5) -- (0,-1);
			}
		}
	}

	\begin{figure}[t]
		\begin{center}
			\begin{tikzpicture}[scale=0.7, transform shape]
				%\draw[help lines] [color=green] ( -5,-5 ) grid ( 5, 5); 
				%\draw (-1,0) -- (1,0); 
				%\draw (0,-1) -- (0,1); 
				\tikzstyle{every node}=[font=\large]
				%Bus1			
				\draw (1,1) pic[thick,scale=1,local bounding box=gen] {Thermal};
				\node[above] at (1,1.6) {P1 $\sim$ P4};
				\draw [very thick] (-1.8,0) -- (1.8,0);	
				\draw [-triangle 45][thick] (-1,0) -- (-1,1.4);
				\node[above] at (-1,1.6) {C1 $\sim$ C4};						
			\end{tikzpicture}
		\end{center}
		\vspace{-4 mm}
		\caption{Diagram for the illustrative example}
		\label{fig:Diagram1}
	\end{figure}
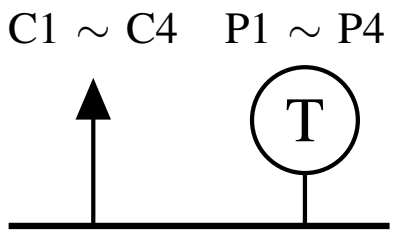

	In this section, we illustrate how RCs impact market equilibrium through a simple example. This example spans six periods corresponding to the  sunset hours, when the net-load increases. 
	
	\subsection{Data}
	
	The market is composed of four producers and four consumers. 
	For the sake of simplicity and clarity, we neglect the transmission system. 
	The four producers are thermal generation units with ramp limits.  
	The system diagram is shown in Fig. \ref{fig:Diagram1}. 
	We consider four consumers representing demands whose utilities are similar.

	Data for producers and consumers are provided in Table \ref{tab:IE_Data}. 
	Columns 2 and 3 include producers' capacities and marginal costs, respectively. 
	Column 4 and 5 provide the ramping limits and the initial production level of each producer. 
	Column 6 represents consumers' utilities for each period, which do not change with time. 
	Table \ref{tab:IE_ConsData} provides consumers' consumption cap in each period. 
	To analyze the impact of RCs clearly, we set a net-load ``jump'' in period 4. 
	The minimum production/consumption levels are set to 0 MW. 
	Moreover, the offer and bid price range is between \$15/MWh and \$0/MWh.

	\begin{table}[t]
		\addtolength{\tabcolsep}{-1.5pt} 
		\centering
		\caption{Producers and Consumers Data}
		\begin{tabular}{c|cccc|c}
			\toprule
			$i$/$j$  & $X_{i,t}^{\max}$ & $c_i$ & $R_i^{\UP}/R_i^{\DN}$ & $x_{i,0}$ &  $u_{j,t}$\\
			\# & (MW) & (\$/MWh) & (MW/h) & (MW) &  (\$/MWh)\\
			&&&&&\\[-2.25ex]\hline &&&&&\\[-1.75ex]
			1 	  & 50    & 0     & 10    & 50       & 13 \\
			2     & 50    & 3     & 10    & 0       & 12 \\
			3     & 50    & 6     & 10    & 0      & 11 \\
			4     & 50    & 9    & 10    & 0     & 2 \\
			\bottomrule
		\end{tabular}%
		\label{tab:IE_Data}%
	\end{table}%
	
	\begin{table}[t]
		\centering
		\caption{Consumers Data: Consumption Cap}
		\begin{tabular}{ccccccc}
			\toprule
			Period (t) & 1    & 2   & 3   & 4  & 5   & 6 \\
			\midrule
			$Q_{j,t}^{\max}$ (MW) & 14    & 15    & 16    & 36    & 35    & 34 \\
			\bottomrule
		\end{tabular}%
		\label{tab:IE_ConsData}%
	\end{table}%
	
	\subsection{Solution: Competitive Equilibria}

	To illustrate the impact of RCs, we first consider the competitive case. 
	Producers/consumers offer/bid at their marginal costs/utilities and true production/consumption caps. 
	We use SW (11a) as the objective function to solve problem (10), with and without RCs. 

		\begin{figure}[t]
		\centering
		\resizebox{0.5\textwidth}{!}{%
		\begin{tikzpicture}
			\pgfplotsset{xmin=0.5, xmax=6.5,
				xtick={1, 2,...,6},
				grid=both, grid style={dash pattern=on 1pt off 1pt},
				%					major grid style={line width=.2pt,draw=gray!50},
				xlabel=Period $t$,
				xlabel near ticks,
				ylabel near ticks,	
				%				label style={font=\small},
				tick label style={font=\small}
			}  
			\begin{axis}[style=semithick, 
				axis y line*=left,
				ymin=-8, ymax=18,
				ylabel=Price (\$/MWh),
				ytick={-10,-8, ...,16}]
				
				\addplot[color=blue, mark=square*, style=very thick]                     
				coordinates
				{
					(1,2) (2,2) (3,2)
					(4,6) (5,6) (6,6)
				};\label{plot_one}					
				
				\addplot [color=red, mark=*, style=very thick]                   
				coordinates
				{
					(1,2) (2,0) (3,-4)
					(4,11) (5,6) (6,6)
				};\label{plot_two}
				
			\end{axis}				
			\begin{axis}[style=semithick, 
				axis x line=none,
				grid = none,
				axis y line*=right,
				ymin=-80, ymax=180,
				ylabel= Quantity (MW),
				label style={font=\small},
				ytick={-100, -80, ...,160},
				legend style={at={(0.94,0.35)}},
				legend cell align={left}]				
				
				\addplot[color=Eggplant, dashed, mark=*, mark options={solid}, style=very thick]                 
				coordinates
				{
					(1,56) (2,60) (3,64)
					(4,144) (5,140) (6,136)
				};
				\addlegendentry{Net-load cap}[\footnotesize]
				\addlegendimage{/pgfplots/refstyle=plot_one}\addlegendentry{Price-NoRCs }[\footnotesize]
				\addlegendimage{/pgfplots/refstyle=plot_two}\addlegendentry{Price-RCs }[\footnotesize]	
			\end{axis}							
		\end{tikzpicture}
		}
		\caption{Competitive Equilibria: Prices}
		\label{fig:C_Price}
	\end{figure}
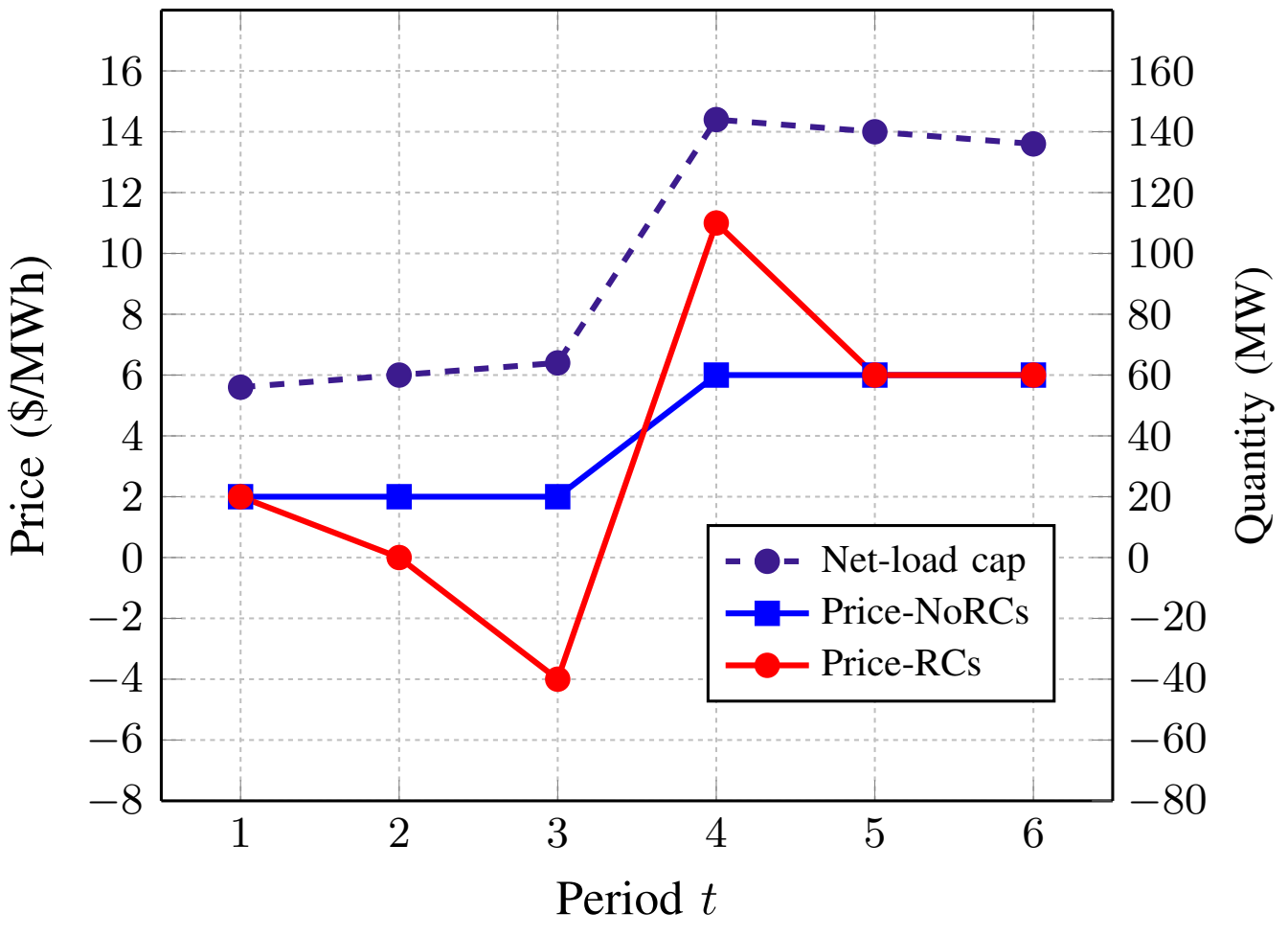

	Fig. \ref{fig:C_Price} shows the price curves for equilibrium A1 (Price-NoRCs) and equilibrium A2 (Price-RCs). 
	The dashed line corresponds to the net-load of the system, which has a ``jump'' in period 4. 
	The price curve of the case with RCs with respect to that of the case without RCs drops down before the net-load cap increase and lifts up in period 4. 
	As the net-load evolves, this price trend is anti-symmetrical during the morning (net-load decrease) and the afternoon (net-load increase) periods.

	From the mathematical viewpoint, some RCs are binding during ramping periods and impact LMPs through constraints (\ref{eq:6b}) and (\ref{eq:6c}). 
	Note that dual variables $\delta_{i,t}^{\UP}$ and $\delta_{i,t}^{\DN}$ link adjacent prices and cause price fluctuations. 
	To help understanding the formation of prices, Table \ref{tab:A2_P2} provides dual variable values for producer 2, which is the marginal producer before the net-load jumps up, 
	We note that when thermal producers proceed to increase their outputs, some dual variables $\delta_{i,t}^{\UP}$ are positive in adjacent periods, which lead to lower/higher price.

		% Table generated by Excel2LaTeX from sheet 'Sheet2'
	\begin{table}[t]
		\centering
		\caption{Equilibrium A2: Results for Producer 2}
		\begin{tabular}{ccccccc}
			\toprule
			Period $t$ &$\lambda_{t}$ &  $c_{2,t}$ & $\delta_{2,t}^{\Xmin}$ & $\delta_{2,t}^{\Xmax}$ & $\delta_{2,t}^{\DN}$   & $\delta_{2,t}^{\UP}$ \\
			\midrule
			3     & -4     & 3     & 0     & 0     & 0     & 4 \\
			4    & 11     & 3     & 0     & 0     & 0     & 11 \\
			5    & 6  		& 3     & 0     & 0     & 0     & 3 \\
			\bottomrule
		\end{tabular}%
		\label{tab:A2_P2}%
	\end{table}%

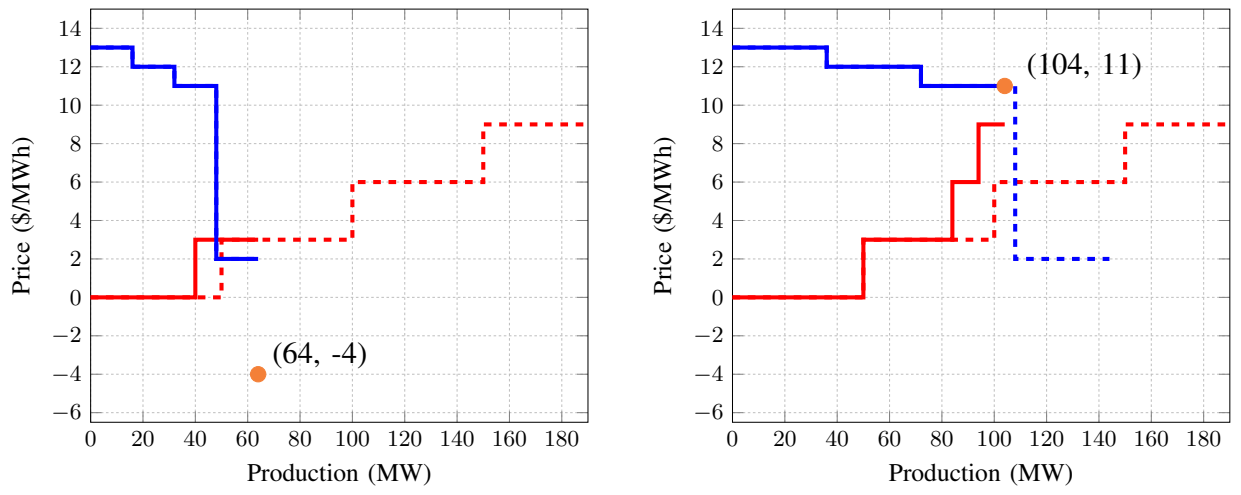
\begin{figure}[t]
		\centering
		\subfigure[Period 3]
		{
			\centering
			\begin{tikzpicture}[scale=0.48]			    
				\pgfplotsset{xmin=0, xmax=190,
					ymin=-6.5, ymax=15,
					ytick={-6, -4, -2, 0, 2, 4, ...,14},
					xtick={0,20,40,...,200},
					grid=both, grid style={dash pattern=on 1pt off 1pt},
					xlabel=Production (MW),
					ylabel=Price (\$/MWh),
					xlabel near ticks,
					ylabel near ticks,	
					label style={font=\normalsize},
					tick label style={font=\small}
				}
				
				\begin{axis}    
					%\addplot[ultra thick, dashed, mark=none,darkgray]  coordinates{ (50,0)     (50,-7) };                  				 							
					\addplot[ultra thick, dashed, mark=none,red]  coordinates{ 	(0,0)     (50,0)  
						(50,3)   (100,3) 
						(100,6)   (150,6)
						(150,9)  (200,9)};	
					%				\addlegendentry{\footnotesize  Demand curve}
					\addplot[ultra thick, dashed, mark=none,blue]  coordinates{ (0,13)     (16,13)  
						(16,12)   (32,12) 
						(32,11)   (48,11)
						(48,2)  (64,2)};	
					%				\addlegendentry{\footnotesize  Demand curve}				  
					\addplot[ultra thick, mark=none,red]  coordinates{ (0,0)     (40,0)  
						(40,3)   (64,3)};	
					\addplot[ultra thick, mark=none,blue]  coordinates{ (0,13)     (16,13)  
						(16,12)   (32,12) 
						(32,11)   (48,11)
						(48,2)	  (64,2)};	
					%				\addlegendentry{\footnotesize  Demand curve}	
					\addplot[only marks, mark=*,mark options={scale=1.5, Orange}] coordinates {(64,-4)};
					%					\addplot[thick, dashed, mark=none,darkgray]  coordinates{ (0,2)     (115,2) };
					%					\addplot[very thick, dashed, mark=none,darkgray]  coordinates{ (100,11)     (100,0) };
					\node[right] at (axis cs: 66,-3) {\large (64, -4)};				
				\end{axis}                         		                       			
			\end{tikzpicture}
		}	
		\subfigure[Period 4]
		{
			\centering
			\begin{tikzpicture}[scale=0.48]			    
				\pgfplotsset{xmin=0, xmax=190,
					ymin=-6.5, ymax=15,
					ytick={-6, -4, -2, 0, 2, 4, ...,15},
					xtick={0,20,40,...,200},
					grid=both, grid style={dash pattern=on 1pt off 1pt},
					xlabel=Production (MW),
					ylabel=Price (\$/MWh),
					xlabel near ticks,
					ylabel near ticks,	
					label style={font=\normalsize},
					tick label style={font=\small}
				}
				
				\begin{axis}
					%\addplot[ultra thick, dashed, mark=none,darkgray]  coordinates{ (100,-7)     (100,11) };               				 							
					\addplot[ultra thick, dashed, mark=none,red]  coordinates{ 	(0,0)     (50,0)  
						(50,3)   (100,3) 
						(100,6)   (150,6)
						(150,9)  (200,9)};	
					%				\addlegendentry{\footnotesize  Demand curve}
					\addplot[ultra thick, dashed, mark=none,blue]  coordinates{ (0,13)     (36,13)  
						(36,12)   (72,12) 
						(72,11)   (108,11)
						(108,2)  (144,2)};	
					%				\addlegendentry{\footnotesize  Demand curve}				  
					\addplot[ultra thick, mark=none,red]  coordinates{ (0,0)     (50,0)  
						(50,3)   (84,3)
						(84,6)   (94,6)	
						(94,9)   (104,9)				};	
					\addplot[ultra thick, mark=none,blue]  coordinates{  (0,13)     (36,13)  
						(36,12)   (72,12)
						(72,11)   (104,11)};	
					%				\addlegendentry{\footnotesize  Demand curve}	
					\addplot[only marks, mark=*,mark options={scale=1.5, Orange}] coordinates {(104,11)};
					%					\addplot[thick, dashed, mark=none,darkgray]  coordinates{ (0,2)     (115,2) };
					%					\addplot[very thick, dashed, mark=none,darkgray]  coordinates{ (100,11)     (100,0) };
					\node[right] at (axis cs: 109,12) {\large (104, 11)};				
				\end{axis}                         		                       			
			\end{tikzpicture}
		}	
		\caption{Equilibrium A2}
		\label{fig:RC1}
	\end{figure}

	From the viewpoint of the market operator, the price fluctuation helps increase production levels, which leads to higher social welfare. 
	Fig. \ref{fig:RC1} shows the offer/bid in periods 3 and 4 using red/blue dashed lines and market clearing results using solid lines for equilibrium A2. 
	We note that the market clearing solutions are not located at the intersection points of offer and bid curves, which is different from the case without RC. 
	To meet increasing net-load in period 4, producer 2 ``over generates'' in period 3. 
	As a result, the price in period 3 is comparatively lower and the market clearing solution locates below the supply curve. 
	Affected by RCs, a ``shortage'' situation happens in period 4. 
	The price in period 4 equals the utility of the marginal consumer. 
	The production and consumption levels of equilibrium A2 are provided in Table \ref{tab:A2_Level}.

	We note that in the competitive case, solving equation (10) is equivalent to solving the KKT conditions of a single market-clearing problem.  
	Dual degeneracy may exist in this problem with specific data set, which results in multiple dual solutions with the same primal solution. 
	This may mislead the equilibrium analysis.

	% Table generated by Excel2LaTeX from sheet 'Illustrative'
	\begin{table}[htbp]
		\centering
		\caption{Equilibrium A2: Production and Consumption Level}
		\resizebox{\columnwidth}{!}{%
		\begin{tabular}{ccccccccccccc}
			\toprule
			$i/j$ & \multicolumn{6}{c}{Production}                 & \multicolumn{6}{c}{Consumption} \\
			\cmidrule(lr){2-7}\cmidrule(lr){8-13} 
			\#    & $t_1$    & $t_2$    & $t_3$    & $t_4$    & $t_5$    & $t_6$    & $t_1$    & $t_2$    & $t_3$    & $t_4$    & $t_5$    & $t_6$     \\
			\midrule
			1     & 50    & 46    & 40    & 50    & 50    & 50    & 14    & 15    & 16    & 36    & 35    & 34 \\
			2     & 4     & 14    & 24    & 34    & 44    & 50    & 14    & 15    & 16    & 36    & 35    & 34 \\
			3     & 0     & 0     & 0     & 10    & 11    & 2     & 14    & 15    & 16    & 32    & 35    & 34 \\
			4     & 0     & 0     & 0     & 10    & 0     & 0     & 12    & 15    & 16    & 0     & 0     & 0 \\
			\bottomrule
		\end{tabular}%
		}
		\label{tab:A2_Level}%
	\end{table}%

	\subsection{Solution: Strategic Equilibria}
	For clarity, we consider only one strategic agent in each
	equilibrium.

	\subsubsection{Strategic Producer} 
	Strategic producers offer to increase the clearing prices.
	We use objective function \eqref{eq:11d} to obtain equilibria with producer $i$ as the strategic producer. 
	Since in equilibrium A2 producers 2 and 3 are marginal producers before or after the net-load jump, we select producer 2 and 3 as strategic producers. 
	
	% Table generated by Excel2LaTeX from sheet 'SingleSP'
	\begin{table}[t]
		\centering
		\caption{Equilibria with Strategic Producer: Price and Offer Price}
		\resizebox{\columnwidth}{!}{%
			\begin{tabular}{lcccccccccccc}
				\toprule
				\multirow{2}[2]{*}{Equilibrium}&  \multicolumn{6}{c}{Prices}                    & \multicolumn{6}{c}{Offer Price} \\
				\cmidrule(lr){2-7}\cmidrule(lr){8-13}   & $t_1$    & $t_2$    & $t_3$    & $t_4$    & $t_5$    & $t_6$   &  $t_1$    & $t_2$    & $t_3$    & $t_4$    & $t_5$    & $t_6$ \\
				\midrule
				A2    & 2     & 0     &-4     & 11    & 6     & 6     & -     & -     & -     & -    & -    & - \\
				B1 (Producer 2)   &  2     & 0     & 0     & 11    & 6     & 6       & 3     & 3     & 7    & 3     & 3    & 3     \\
				B2 (Producer 3)   & 2     & 0     & -7    & 11    & 9     & 9     & 6     & 6     & 6    & 6     & 9    & 9 \\
				\bottomrule
			\end{tabular}%
		}
		\label{tab:SP_PO}%
	\end{table}%

	The price and offer results for equilibria with strategic producers are provided in Table \ref{tab:SP_PO}. 
	RCs are considered in each equilibrium for comparison. 
	The first row provides the prices of equilibrium A2. 
	The second row provides the prices for equilibrium B1 and the offer prices for producer 2. 
	In equilibrium B1, producer 2 is the strategic producer and others agents are competitive. 
	We note that producer 2 offer high in period 3, which lifts up the market clearing price in that period.
	The third row provides the prices for equilibrium B2 and producer 3's offer prices. 
	In equilibrium B2, producer 3 is the strategic agent. 
	Producer 3 offers high in periods 5 and 6, which lead to higher price in these periods but lower price in period 3. 

	We note that the strategic producer that produces before and after the net-load jump in the competitive case offer prices strategically to reduce the price fluctuation. 
	This behavior reduces its intertemporal opportunity cost. 
	On the other hand, the strategic producer that only produces after the net-load jump in the competitive case offers prices strategically to raise the price after the net-load jump. 
	This behavior brings more distortion to the clearing price in periods before the net-load jump and leads to higher intertemporal opportunity costs for other producers.

	\subsubsection{Strategic Consumers}
	Strategic consumers bid to reduce prices.
	We use objective function (\ref{eq:11e}) to obtain the equilibria with consumer $j$ as the strategic consumer.  
	Since in equilibrium A2 consumers 3 and 4 are marginal consumers after or before the net-load jump, we select consumers 3 and 4 as strategic consumers. The prices and bids results are shown in Table \ref{tab:SC_PB}. 
	The first row provides the price of equilibrium B3 and the bid prices for consumer 3. 
	Consumer 3 bids lower after the net-load jump up, which lead to lower prices in periods 4-6. 
	Consequentially, the price in period 3 raise up, which reduces intertemporal opportunity costs for producers.
	The second row provides the clearing prices for equilibrium B4 and bid prices for consumer 4.
	Consumer 4 bids at its utility in all periods and doesn't alter the market clearing price.

		% Table generated by Excel2LaTeX from sheet 'SingleSP'
	\begin{table}[t]
		\centering
		\caption{Equilibria with Strategic Consumer: Price and Bid Price}
		\resizebox{\columnwidth}{!}{%
			\begin{tabular}{lcccccccccccc}
				\toprule
				\multirow{2}[2]{*}{Equilibrium}&  \multicolumn{6}{c}{Prices}                    & \multicolumn{6}{c}{Bid Price} \\
				\cmidrule(lr){2-7}\cmidrule(lr){8-13}   & $t_1$    & $t_2$    & $t_3$    & $t_4$    & $t_5$    & $t_6$   &  $t_1$    & $t_2$    & $t_3$    & $t_4$    & $t_5$    & $t_6$ \\
				\midrule
				
				B3 (Consumer 3)   &  2     & 0     & -2     & 9    & 6     & 3      & 11     & 11    & 11    & 9    & 6   & 3     \\
				B4 (Consumer 4)   & 2     & 0     &-4     & 11    & 6     & 6     & 2     & 2     & 2     & 2   & 2    & 2 \\
				\bottomrule
			\end{tabular}%
		}
		\label{tab:SC_PB}%
	\end{table}%
	
	Comparing the competitive equilibrium A2 and the strategic equilibria B1-B4, we observe that the impact of RCs on market equilibrium is time coupled. 
	Strategic behaviors of agents in certain periods may affect the market-clearing outcomes in adjacent periods, even if the agent does not offer/bid strategically in those periods.  
	The strategic behavior before or after a net-load jump may lead to a ``seesaw'' behavior on the market-clearing price. 
	Strategically increasing prices after a net-load jump may lead to lower clearing prices before the net-load jump and higher intertemporal opportunity costs for producers, and vice versa.

	\section{Case Study}
	\label{CaseStudy}
	
	In this section, we provide results from a case study based on the IEEE Reliability Test System (RTS) \cite{Christos(2016)}. 
	We consider only thermal generating units in the system and simulate 8 time periods spanning the zone with high net-demand variability.
	
	\subsection{Data Set}
	\label{Data Set}

	Table \ref{tab:case_pc} provides producer and consumer data, where $c_i$, $ x_{i,0}$, and $X_i^{\max}$ correspond to the cost, initial production and maximum production level (capacity) of producer $i$ and $u_j$ and $Q_j^{\max}$ correspond to base values for utility and maximum consumption level for consumer $j$, obtained from consumption data in \cite{Christos(2016)}.

	\begin{table*}[t!]
		\centering
		\caption{Producer and Consumer Data}
		\resizebox{2\columnwidth}{!}{%
			\begin{tabular}{|c|c|c|c|c|c|c|c|c|c|c|c|c|c|c|c|c|c|}
				\hline
				Producer   & 1     & 2     & 3     & 4     & 5     & 6     & 7     & 8     & 9     & 10    & 11    & 12    &      \multicolumn{5}{c|}{\multirow{4}[8]{*}{}}   \\
				\cline{1-13}
				$X_i^{\max}$ (MW) & 152   & 152   & 350   & 591   & 60    & 155   & 155   & 400   & 400   & 300   & 310   & 350   &  \multicolumn{5}{c|}{} \\
				\cline{1-13}
				$c_i$ (\$/MWh)  & 13.32 & 15.32 & 19.7  & 20.93 & 26.11 & 6.52  & 8.52  & 5.02  & 3.47  & 0     & 11.52 & 12.49 &      \multicolumn{5}{c|}{} \\
				\cline{1-13}
				$x_{i,0}$ (MW) & 0     & 0     & 0     & 0     & 0     & 0     & 0     & 0     & 246   & 300   & 0     & 0     &      \multicolumn{5}{c|}{} \\
				\hline
				\hline
				Consumer  & 1     & 2     & 3     & 4     & 5     & 6     & 7     & 8     & 9     & 10    & 11    & 12    & 13    & 14    & 15    & 16    & 17\\
				\hline
				$Q_j^{\max}$ (MW) & 101   & 90    & 167   & 69    & 66    & 127   & 117   & 159   & 162   & 180   & 246   & 180   & 294   & 93    & 310   & 170   & 119 \\
				\hline
				$u_j$ (\$/MWh)  & 50    & 49    & 48    & 47    & 46    & 45    & 44    & 43    & 42    & 41    & 35    & 30    & 25    & 20    & 15    & 7     & 2 \\
				\hline
			\end{tabular}%
		}
		\label{tab:case_pc}%
	\end{table*}%

\begin{figure}[t!]
		\centering
				\subfigure[Consumption Ratio ($r_{j,t}^{\rm NL}$)]
		{
			\centering
			\begin{tikzpicture}[scale=0.48]			    
				\pgfplotsset{xmin=0, xmax=24,
					ymin=0, ymax=1.2,
					ytick={0.1, 0.2, ..., 1.2},
					xtick={0,2,4,...,24},
					grid=both, grid style={dash pattern=on 1pt off 1pt},
					xlabel= Period ($t$),
					ylabel= Consume Ratio (\%),
					xlabel near ticks,
					ylabel near ticks,	
					label style={font=\normalsize},
					tick label style={font=\small}
				}				
				\begin{axis}                  				 					
					\addplot[ultra thick, mark=none,red]  coordinates{ 						
						%						(0, 0.72)
						(1, 0.71)
						(2, 0.68)
						(3, 0.69)
						(4, 0.70)
						(5, 0.73)
						(6, 0.78)
						(7, 0.84)
						(8, 0.60)
						(9, 0.36)
						(10, 0.27)
						(11, 0.24)
						(12, 0.19)
						(13, 0.18)
						(14, 0.16)
						(15, 0.17)
						(16, 0.22)
						(17, 0.32)
						(18, 0.62)
						(19, 0.89)
						(20, 0.91)
						(21, 0.86)
						(22, 0.84)
						(23, 0.80)
					};	
					%				\addlegendentry{\footnotesize  Demand curve}													

					\fill[red, opacity=0.15, draw=none] (axis cs:5,0) rectangle (axis cs:12,1.2);
					
				\end{axis}                         		                       		
			\end{tikzpicture}
			\label{fig:NetLoad}
		}
		\subfigure[Utility Ratio ($r_{j,t}^{\rm UT}$)]
		{
			\centering
			\begin{tikzpicture}[scale=0.48]			    
				\pgfplotsset{xmin=0, xmax=24,
					ymin=0, ymax=1.2,
					ytick={0.1, 0.2, ..., 1.2},
					xtick={0,2,4,...,24},
					grid=both, grid style={dash pattern=on 1pt off 1pt},
					xlabel= Period ($t$),
					ylabel= Consume Ratio (\%),
					xlabel near ticks,
					ylabel near ticks,	
					label style={font=\normalsize},
					tick label style={font=\small}
				}				
				\begin{axis}                  				 			
					\addplot[ultra thick, mark=none,blue]  coordinates{ 
%						(0, 0.83)
						(1, 0.80)
						(2, 0.78)
						(3, 0.77)
						(4, 0.77)
						(5, 0.80)
						(6, 0.86)
						(7, 0.93)
						(8, 0.94)
						(9, 0.90)
						(10, 0.85)
						(11, 0.80)
						(12, 0.76)
						(13, 0.74)
						(14, 0.74)
						(15, 0.75)
						(16, 0.79)
						(17, 0.84)
						(18, 0.92)
						(19, 0.97)
						(20, 1.00)
						(21, 0.97)
						(22, 0.95)
						(23, 0.89)
					};	
					%				\addlegendentry{\footnotesize  Demand curve}		
					
					\fill[blue, opacity=0.15, draw=none] (axis cs:5,0) rectangle (axis cs:12,1.2);		
				\end{axis}                         		                       		
			\end{tikzpicture}
			\label{fig:ConsRate}
		}

		\caption{Production and Consumption Ratio }
		\label{fig:Orate}
		\end{figure}
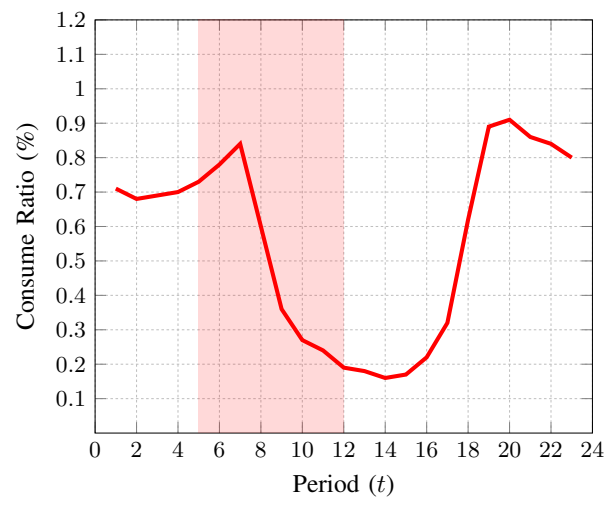
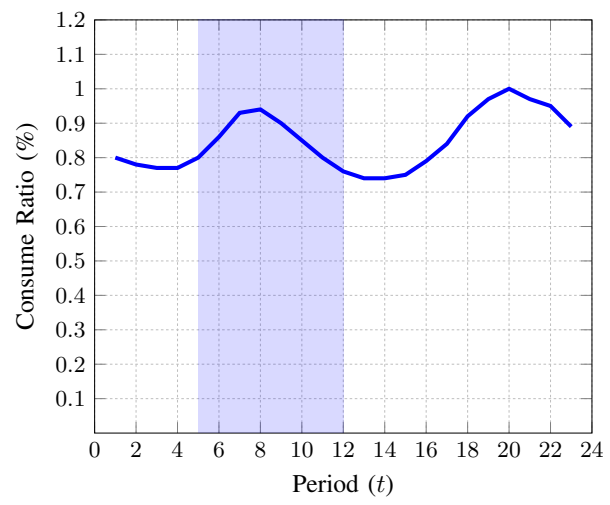

	Focusing on the impact of ramping constraints, we only consider thermal generating units.
	On the other hand, we consider that the consumers are cognizant of the shape of the net demand and seek to take advantage of sharp ramping zones in that demand. 
	The consumption cap for each period is computed as $Q_{j}^{\max}\cdot r_{j,t}^{\rm NL}$, where $r_{j,t}^{\rm NL}$ is a consumption ratio.
	The utility for each periods is computed as $u_{j} \cdot r_{j,t}^{\rm UT}$, where $r_{j,t}^{\rm UT}$ is a utility ratio. 
	The consumption ratio $r_{j,t}^{\rm NL}$ and utility ratio $r_{j,t}^{\rm UT}$ are derived from CAISO net-load and demand data on April $1^{\rm st}$, 2022 \cite{CAISO, GridStatus}.
	We divide the net-load/demand value by the maximum demand on that day and obtain consumption/utility ratios.
	The $r_{j,t}^{\rm NL}$ and $r_{j,t}^{\rm UT}$ curves for consumer $j$ are shown in Fig. \ref{fig:NetLoad} and Fig. \ref{fig:ConsRate}, respectively.
	For the simulation, we focus on the sunrise interval, from hour 5 to hour 12, spanning 8 time periods.
	Additionally, the hourly ramping-up/down limits for thermal power producers are set to 30\% respecting to their respective capacities.
	The offer and bid prices range between \$50/MWh and \$0/MWh.

	\subsection{Results}
	\label{Results}

	\begin{table}[t!]
		\centering
		\caption{Case Study: Market Outcomes}
		\addtolength{\tabcolsep}{-2.2pt}   
		\begin{tabular}{lcccccc}
			\toprule
			\multirow{2}{*}{Equilibrium} & \multicolumn{1}{c}{Max} & SW    & TP    & TS    & Production & Avg. Price \\
			&\multicolumn{1}{c}{OF} & (\$)  & (\$)  & (\$)  & (MW)  & (\$/MWh) \\
			\midrule
			C1(NoRCs) & SW    & 249,093 & 51,500 & 197,593 & 9,671 & 7.99  \\
			C2(RCs) & SW     & 247,072 & 48,686 & 198,386 & 9,493 & 7.70  \\
			\midrule
			D1(NoRCs) & SW    & 249,093 & 63,867 & 185,226 & 9,671 & 9.52  \\
			D2(RCs) & SW    & 247,072 & 55,803 & 191,269 & 9,493 & 8.66  \\
			\midrule
			D3(NoRCs) & TP   & 191,144 & 167,422 & 23,722 & 5,238 & 32.60  \\
			D4(RCs) & TP    & 193,277 & 154,051 & 39,226 & 5,704 & 30.27  \\
			\midrule
			D5(NoRCs) & TS   & 238,090 & 22,863 & 215,227 & 7,867 & 5.44  \\
			D6(RCs) & TS     & 233,511 & 17,274 & 216,237 & 8,153 & 5.01  \\
			\bottomrule
		\end{tabular}%
		\label{tab:Case_MO}%
	\end{table}%

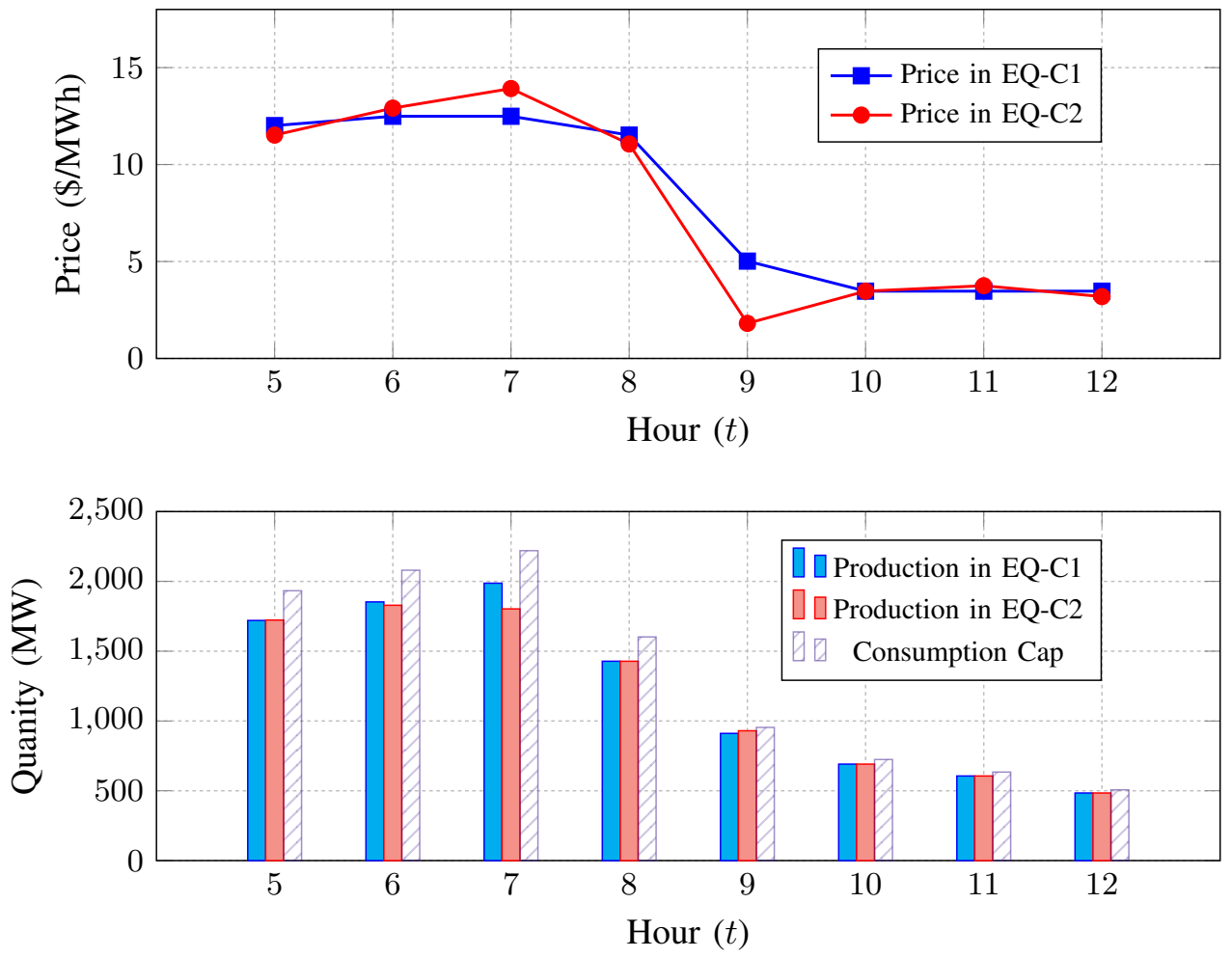
\begin{figure}[t!]
		\centering
	
		\begin{tikzpicture}[scale=0.7]
			\begin{axis}[name=plot1,
				height=5cm,width=12cm,
				xmin=4, xmax=13,
				ymin=0, ymax=18,
				ytick={0, 5, ..., 20},
				xtick={5,6,...,10,11,12},
				grid=both, grid style={dash pattern=on 1pt off 1pt},
				xlabel= Hour ($t$),
				ylabel= Price (\$/MWh),
				xlabel near ticks,
				ylabel near ticks,	
				label style={font=\normalsize},
				tick label style={font=\small},
				legend style={at={(axis cs:12,16)},font=\footnotesize},
				]
				
				\addplot[thick, mark=square*,blue]  coordinates{ 
					%						(0, 0.83)
					(5, 12.01)
					(6, 12.49)
					(7, 12.49)
					(8, 11.52)
					(9, 5.02)
					(10, 3.47)
					(11, 3.47)
					(12, 3.47)
				};	
				%				\addlegendentry{\footnotesize  Demand curve}		
				\addplot[thick, mark=*,red]  coordinates{ 
					%						(0, 0.83)
					(5, 11.52)
					(6, 12.91)
					(7, 13.92)
					(8, 11.06)
					(9, 1.81)
					(10, 3.47)
					(11, 3.75)
					(12, 3.19)
				};	
				%				\addlegendentry{\footnotesize  Demand curve}		
				
				\legend {Price in EQ-C1, Price in EQ-C2};
				
			\end{axis}    
			
			\begin{axis}[name=plot3,at={($(plot1.south)-(0,1.5cm)$)},anchor=north,
				height=5cm,width=12cm,
				xmin=4,xmax=13,
				ymin=0,ymax=2500,
				xtick={5,6,7,...,11,12},	
				ytick={0,500,...,2500},	
				grid=both, grid style={dash pattern=on 1pt off 1pt},
				xlabel= Hour ($t$),
				ylabel= Quanity (MW),
				xlabel near ticks,
				ylabel near ticks,	
				label style={font=\normalsize},
				tick label style={font=\small},
				legend style={at={(axis cs:12,2300)},font=\footnotesize},
				ybar,
				bar width=5pt,
				ybar=0pt,
				xtick align=inside
				]
				\addplot [draw=blue,fill=cyan] coordinates {
					(5,1720) 
					(6,1853) 
					(7,1987) 
					(8,1428) 
					(9,911) 
					(10,691) 
					(11,606) 
					(12,484) 
				};
				\addplot [draw=red,fill=Salmon]coordinates {
					(5,1723) 
					(6,1829) 
					(7,1803) 
					(8,1428) 
					(9,930) 
					(10,691) 
					(11,606) 
					(12,484) 
				};
				
				\addplot [draw=Eggplant, pattern =north east lines, pattern color = Eggplant, opacity=0.5]coordinates {
					(5,1933) 
					(6,2080) 
					(7,2220) 
					(8,1602) 
					(9,954) 
					(10,724) 
					(11,634) 
					(12,507) 
				};

				\legend {Production in EQ-C1, Production in EQ-C2, Consumption Cap};				
			\end{axis}   
			
		\end{tikzpicture}
	\caption{Equilibria C1 and C2: Price and Production Level}
	\label{fig:EQ-C}
	\end{figure}

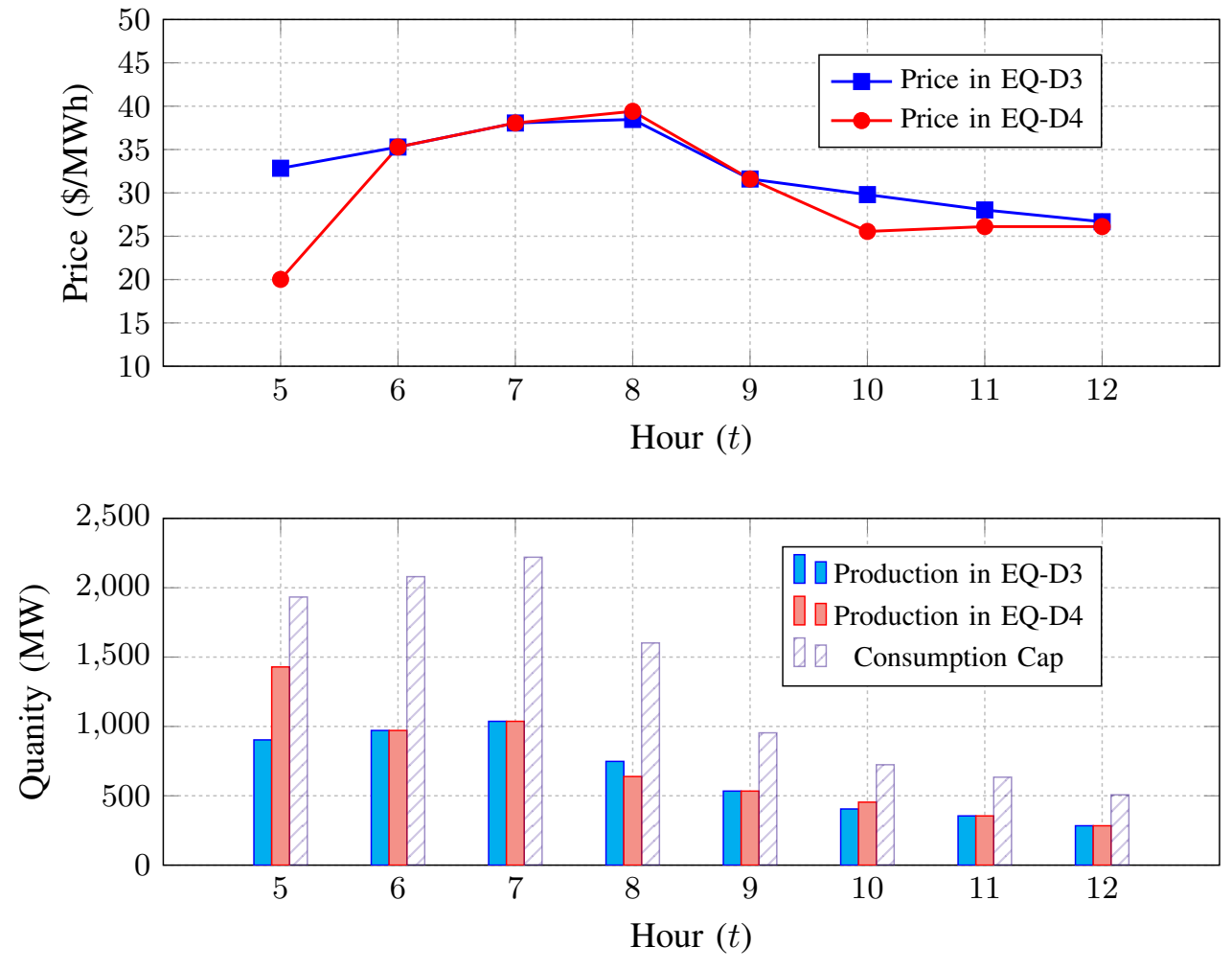
\begin{figure}[t!]
		\centering
		
		\begin{tikzpicture}[scale=0.7]
			\begin{axis}[name=plot1,
				height=5cm,width=12cm,
				xmin=4, xmax=13,
				ymin=10, ymax=50,
				ytick={0, 5, ..., 50},
				xtick={5,6,...,10,11,12},
				grid=both, grid style={dash pattern=on 1pt off 1pt},
				xlabel= Hour ($t$),
				ylabel= Price (\$/MWh),
				xlabel near ticks,
				ylabel near ticks,	
				label style={font=\normalsize},
				tick label style={font=\small},
				legend style={at={(axis cs:12,46)},font=\footnotesize},
				]
				
				\addplot[thick, mark=square*,blue]  coordinates{ 
					%						(0, 0.83)
					(5, 32.84)
					(6, 35.29)
					(7, 38.06)
					(8, 38.47)
					(9, 31.60)
					(10, 29.80)
					(11, 28.03)
					(12, 26.67)
				};	
				%				\addlegendentry{\footnotesize  Demand curve}		
				\addplot[thick, mark=*,red]  coordinates{ 
					%						(0, 0.83)
					(5, 20.02)
					(6, 35.29)
					(7, 38.06)
					(8, 39.41)
					(9, 31.60)
					(10, 25.55)
					(11, 26.11)
					(12, 26.11)
				};	
				%				\addlegendentry{\footnotesize  Demand curve}		
				
				\legend {Price in EQ-D3, Price in EQ-D4};
				
			\end{axis}    
			
			\begin{axis}[name=plot3,at={($(plot1.south)-(0,1.5cm)$)},anchor=north,
				height=5cm,width=12cm,
				xmin=4,xmax=13,
				ymin=0,ymax=2500,
				xtick={5,6,7,...,11,12},	
				ytick={0,500,...,2500},	
				grid=both, grid style={dash pattern=on 1pt off 1pt},
				xlabel= Hour ($t$),
				ylabel= Quanity (MW),
				xlabel near ticks,
				ylabel near ticks,	
				label style={font=\normalsize},
				tick label style={font=\small},
				legend style={at={(axis cs:12,2300)},font=\footnotesize},
				ybar,
				bar width=5pt,
				ybar=0pt,
				xtick align=inside
				]
				\addplot [draw=blue,fill=cyan] coordinates {
					(5,903) 
					(6,971) 
					(7,1037) 
					(8,748) 
					(9,534) 
					(10,405) 
					(11,355) 
					(12,284) 
				};
				\addplot [draw=red,fill=Salmon]coordinates {
					(5,1429) 
					(6,971) 
					(7,1037) 
					(8,639) 
					(9,534) 
					(10,454) 
					(11,355) 
					(12,284) 
				};
				
				\addplot [draw=Eggplant, pattern =north east lines, pattern color = Eggplant, opacity=0.5]coordinates {
					(5,1933) 
					(6,2080) 
					(7,2220) 
					(8,1602) 
					(9,954) 
					(10,724) 
					(11,634) 
					(12,507) 
				};

				\legend {Production in EQ-D3, Production in EQ-D4, Consumption Cap};				
			\end{axis}   
			
		\end{tikzpicture}
		\caption{Equilibria D3 and D4: Price and Production Level}
		\label{fig:EQ-TP}
	\end{figure}

\begin{figure}[t!]
		\centering
		
		\begin{tikzpicture}[scale=0.7]
			\begin{axis}[name=plot1,
				height=5cm,width=12cm,
				xmin=4, xmax=13,
				ymin=0, ymax=18,
				ytick={0, 5, ..., 20},
				xtick={5,6,7,...,11,12},
				grid=both, grid style={dash pattern=on 1pt off 1pt},
				xlabel= Hour ($t$),
				ylabel= Price (\$/MWh),
				xlabel near ticks,
				ylabel near ticks,	
				label style={font=\normalsize},
				tick label style={font=\small},
				legend style={at={(axis cs:12,16)},font=\footnotesize},
				]
				
				\addplot[thick, mark=square*,blue]  coordinates{ 
					%						(0, 0.83)
					(5, 6.52)
					(6, 6.52)
					(7, 8.52)
					(8, 6.52)
					(9, 5.02)
					(10, 3.47)
					(11, 3.47)
					(12, 3.47)
				};	
				%				\addlegendentry{\footnotesize  Demand curve}		
				\addplot[thick, mark=*,red]  coordinates{ 
					%						(0, 0.83)
					(5, 5.02)
					(6, 6.52)
					(7, 6.52)
					(8, 11.52)
					(9, 1.52)
					(10, 2.02)
					(11, 5.02)
					(12, 1.92)
				};	
				%				\addlegendentry{\footnotesize  Demand curve}		
				
				\legend {Price in EQ-D5, Price in EQ-D6};
				
			\end{axis}    
			
			\begin{axis}[name=plot3,at={($(plot1.south)-(0,1.5cm)$)},anchor=north,
				height=5cm,width=12cm,
				xmin=4,xmax=13,
				ymin=0,ymax=2500,
				xtick={5,6,7,...,11,12},	
				ytick={0,500,...,2500},	
				grid=both, grid style={dash pattern=on 1pt off 1pt},
				xlabel= Hour ($t$),
				ylabel= Quanity (MW),
				xlabel near ticks,
				ylabel near ticks,	
				label style={font=\normalsize},
				tick label style={font=\small},
				legend style={at={(axis cs:12,2300)},font=\footnotesize},
				ybar,
				bar width=5pt,
				ybar=0pt,
				xtick align=inside
				]
				\addplot [draw=blue,fill=cyan] coordinates {
					(5,1255) 
					(6,1255) 
					(7,1410) 
					(8,1255) 
					(9,911) 
					(10,691) 
					(11,606) 
					(12,484) 
				};
				\addplot [draw=red,fill=Salmon]coordinates {
					(5,1569) 
					(6,1441) 
					(7,1286) 
					(8,1333) 
					(9,907) 
					(10,620) 
					(11,606) 
					(12,392) 
				};
				
				\addplot [draw=Eggplant, pattern =north east lines, pattern color = Eggplant, opacity=0.5]coordinates {
					(5,1933) 
					(6,2080) 
					(7,2220) 
					(8,1602) 
					(9,954) 
					(10,724) 
					(11,634) 
					(12,507) 
				};

				\legend {Production in EQ-D5, Production in EQ-D6, Consumption Cap};				
			\end{axis}   
			
		\end{tikzpicture}
		\caption{Equilibria D5 and D6: Price and Production Level}
		\label{fig:EQ-TS}
	\end{figure}
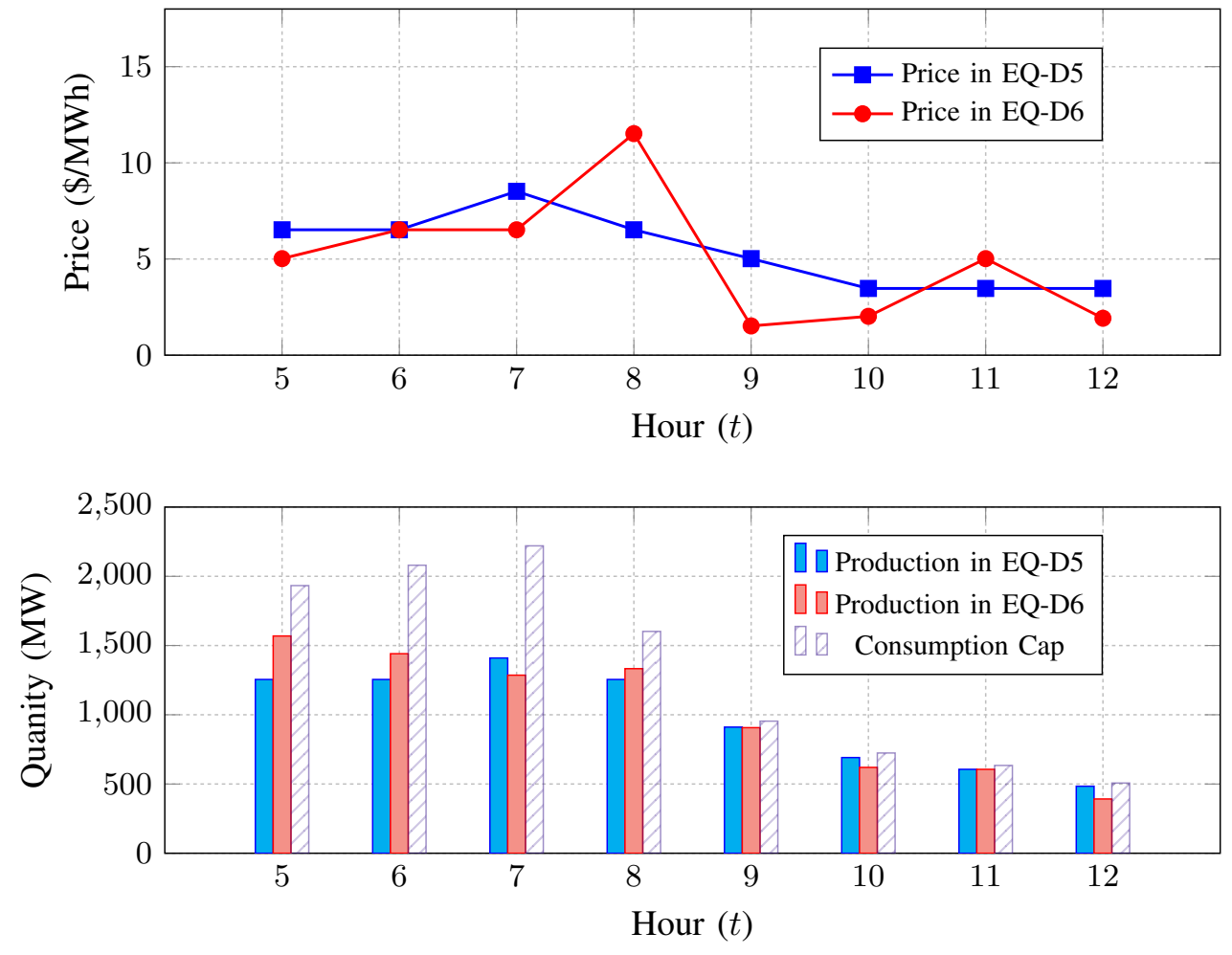

	Market outcomes are provided in Table \ref{tab:Case_MO}.
	Additionally, Fig. \ref{fig:EQ-C}, \ref{fig:EQ-TP}, and \ref{fig:EQ-TS} provide price curves and production bar charts for different equilibria.
    To show the effect of RCs, we use blue/red lines and bars to represent prices and production outcomes with/without RCs.
	We use purple bars with diagonal stripes to represent consumption cap for reference.

	\subsubsection{Competitive Equilibria}
	Equilibria C1 and C2 correspond to competitive equilibria without and with RCs.
	The market outcomes are reported in Table \ref{tab:Case_MO}, where equilibrium C1 and C2 correspond to the case without and with RCs, respectively. 
	We note that the SW and average price in equilibrium C2 are lower than those in equilibrium C1. 
	Regarding the 8-period sunrise simulation, we note that RCs lead to lower TP and production levels and higher TS. 
	As compared with the price curve without RCs, the price curve with RCs presents a lift-and-drop trend, which is anti-symmetrical to the trend during the sunset horizon, as shown in Fig. \ref{fig:C_Price}.
	To meet the ramping need during the sunrise period, lower cost thermal producers reduce generation and higher cost units increase generation, rendering higher price during hours 7-8.
	On the contrary, prices are lower during hours 8-9, stimulating thermal producers to decrease production.
	Comparing with the production profile without RCs, the production profile with RCs is lower during hours 6-7 limited by ramping capabilities, though the consumption increases.

	\subsubsection{Strategic Equilibira}

	If the OF is set to SW, we obtain equilibria D1 (NoRCs) and D2 (RCs).
	We note that the SW and production levels in equilibria D1 and D2 are the same as those in equilibria C1 and C2, respectively. 
	Since multiple equilibria exist \cite{Liu(2021)}, the market clearing price in a strategic equilibrium can be higher or lower than that in the corresponding competitive equilibrium.

	If the objective function of problem (11) is set to TP, we obtain equilibria D3 and D4 with the highest prices. 
	The price curves trend for equilibria D3 and D4 are similar to those of competitive equilibria C1 and C2, respectively.
	At hour 7, the price and production level for equilibria D3 and D4 are the same, while the total profits are different (\$36,374 in D3 and \$35,870 in D4).
	This is because, with ramping constraints, a high cost unit produces more in equilibrium D4.
	Additionally, we note that the TP and the average of equilibrium D4 are lower than those of equilibrium D3, while the production level, SW and TS are higher.
	RCs limit the flexibility of producers to change production in adjacent periods, thus reducing the market power of strategic producers.

	If the objective function is set to TS, we obtain equilibria D5 and D6, in which the prices are comparatively lower. 
	The ``seesaw'' characteristic is observed on the price curve of equilibrium D6.
	The trend of the price curve for equilibrium D6 (RCs) is different from that of production level because of the impact of ramping limits.
	We note that the TS and production level of equilibrium D6 are higher than those in equilibrium D5.
	With RCs, the strategic consumers has higher TS as compared to the case without RCs.
	This indicates that RCs increase the price-manipulation capability of strategic consumers.

	\subsection{Computational Considerations}
	
	The case study results are obtained using the global solver BARON \cite{BARON(2005)} under GAMS on the NEOS server \cite{Czyzyk(1998)}. To reduce the computational burden and obtain high quality results, we implement a two-step solution strategy:
	\begin{enumerate}
		\item We solve single-period problems for each period sequentially and get an initial solution. 
		\item We specify ranges for $\lambda_{t}$ for all periods based on upper/lower bound of offer/bid prices. Then, using the initial solution as starting point, we solve the problem for all periods with a time limit. 
	\end{enumerate}

	If the objective function is SW, which is linear, the competitive case requires much less computational time than the strategic one. 
	On the other hand, the strategic case brings more price flexibility, which leads to multiple equilibria and high computational burden. 
	Using nonlinear objective functions requires significantly higher computational time than using linear objective functions.
	Additionally, RCs may lead to additional multiple equilibria with different production allocations, which brings extra difficulty to the solver.

	\section{Conclusion}
	\label{Conclusion}
	This paper uses a multi-period equilibrium model to study the impact of RCs on electricity market outcomes, particularly clearing prices. 
	In the illustrative example, the drop-and-lift behavior of the price curve during ramping periods is carefully investigated from a mathematical viewpoint and from simulation results.
	We note that RCs generally lead to multiple equilibria in the competitive case. 
	On the other hand, the case study shows that RCs reduce market power of strategic producers, while increase price-manipulation capability for strategic consumers.
	Regarding market outcomes, RCs impact on profits and surpluses is highly dependent on the system data set.

	\begin{IEEEbiographynophoto}
		%		[{\includegraphics[clip,width=1in,height=1.25in,keepaspectratio]{./Figure/liu.pdf}}]
		{Xuan Liu}
		(M'$24$) received the B.Eng.\ degree from the Department of Electrical Engineering, Tsinghua University, Beijing, China, in 2018, and the M.S.\ and the Ph.D. degrees from the Department of Electrical and Computer Engineering, The Ohio State University, Columbus, Ohio, US, in 2019 and 2024, respectively.
		\\\indent 
		He is currently an Engineer with Beijing Power Exchange Center, Beijing, China.	His research interests include the electricity market, power system operations and planning, optimization theory and its applications.
	\end{IEEEbiographynophoto}

	\begin{IEEEbiographynophoto}
		%		[{\includegraphics[clip,width=1in,height=1.25in,keepaspectratio]{./Figure/conejo.pdf}}]
		{Antonio J. Conejo}
		(F'$04$) received the M.S.\ degree from the Massachusetts Institute of Technology, Cambridge, MA, in~$1987$, and the Ph.D.\ degree from the Royal Institute of Technology, Stockholm, Sweden, in~$1990$.\\
		\indent He is currently a Professor in the Department of Integrated Systems Engineering and the Department of Electrical and Computer Engineering, The Ohio State University, Columbus, Ohio. His research interests include control, operations, planning, economics and regulation of electric energy systems, as well as statistics and optimization theory and its applications.
	\end{IEEEbiographynophoto}

	\vfill

	%\balance
\end{document}